\documentclass[11pt]{amsart}
\newtheorem{Theorem}{Theorem}[section]

\newtheorem{Assumption 2}[Theorem]{Assumption 2}

\numberwithin{equation}{section}

\usepackage{float}
\usepackage{amssymb}
\usepackage[dvips]{graphicx}
\usepackage{epsfig}
\usepackage{caption}
\usepackage{yhmath}

\usepackage{stmaryrd}
\usepackage{ulem}
\usepackage{yfonts}
\usepackage{amsbsy}
\usepackage{pstricks}
\usepackage{graphicx}
\usepackage{pst-xkey}
\usepackage{pst-all,pst-gr3d}
\usepackage[dvips,arrow,matrix]{xy}

\def\F{{\mathcal F}}
\def\DD{{\mathbb D}}

\def\KK{{\mathbb K}}
\def\NN{{\mathbb N}}
\def\MM{{\mathbb M}}
\def\P{{\mathbb{P}}}
\def\RR{{\mathbb{R}}}
\def\ZZ{{\mathbb Z}}

\begin{document}
\address{Universit\'e de Bordeaux, Institut de Math\'ematiques, UMR CNRS 5251, F-33405 Talence Cedex}
\email{agnes.bachelot@math.u-bordeaux1.fr}

\title{Wave Computation on the Hyperbolic Double Doughnut}

\author{Agn\`{e}s BACHELOT-MOTET}

\begin{abstract}
We compute the waves propagating on the compact surface of constant
negative curvature and genus 2. We adopt a variational approach using finite elements. We
have to implement the action of the fuchsian group by suitable
boundary conditions of periodic type. A spectral analysis of the wave
allows to compute the spectrum and the eigenfunctions of the
Laplace-Beltrami operator. We test the exponential decay due to a
localized dumping and the ergodicity of the geodesic flow.
\end{abstract}
\maketitle

\pagestyle{myheadings}
\markboth{\centerline{\sc Agn\`{e}s Bachelot-Motet}}{\centerline{\sc
 Wave Computation on the Hyperbolic Double Doughnut}}


\section{Introduction.}

The Hyperbolic Double Doughnut $\mathbf{K}$ is the compact manifold of negative
constant curvature with two
holes. We can define it by the quotient of the hyperbolic Poincar\'e
disc $\mathbf{D}$, by some Fuchsian group $\Gamma$. Alternatively, we
can construct it as the quotient of the so-called Dirichlet polygon,
or fundamental domain $\mathcal{F}\subset\mathbf{D}$ by a suitable
relation of equivalence $\sim$ :
$$
\mathbf{K}=\mathbf{D}/\Gamma=\mathcal{F}/\sim.
$$
 This beautiful object has many fascinating properties as regards the
classical and quantum chaos (classical references are \cite{Voroz}, \cite{gutz}). Several
important 
computational investigations of the spectrum were performed by using a stationnary
method by R. Aurich and F. Steiner \cite{Aurich1}.  Moreover
there has been much recent interest for  the cosmological models with
non trivial topology (a seminal work is the famous ``Cosmic Topology''
by M. Lachi\`{e}ze-Rey and J-P. Luminet \cite{Luminet}). In this
context, $\mathbf{K}$ has been studied as a paradigm in
\cite{Cornish}, where a schema
based on the 
finite differences on a euclidean grid was used to solve the D'Alembertian.
In this paper we compute the solutions of the  wave equations in the
time domain, by using a variational method and a discretization with
finite elements on very fine meshes. The domain of calculus is the
Dirichlet polygon, therefore the initial Cauchy problem on
the manifold without booundary $\mathbf{K}$, becomes a mixed problem
on  $\mathcal{F}$ and the action of the Fuchsian group is expressed as
boundary conditions on $\partial\mathcal{F}$, analogous to periodic
conditions. These boundary constraints are implemented in the choice of
the basis of finite elements. By this way we obtain very accurate
results on the transient waves. We test these results by performing a
Fourier analysis of the transients waves that allows to find the first
eigenvalues of the Laplace-Beltrami operator $\Delta_{\mathbf{K}}$ on $\mathbf{K}$. We
compute also the solutions of the damped wave equation 
$$
\partial_t^2\psi-\Delta_{\mathbf{K}}\psi+a\partial_t\psi=0.
$$
When $0\leq a\in L^{\infty}(\mathbf{K})$ and $a>0$ on
$\partial\mathcal{F}$, the ergodicity of the
geodesic flow assures that the geometric control condition of Rauch
and Taylor \cite{Rauch} is satisfied. Our numerical experiments agree with their
theoretical results, stating that the energy decays exponentially.

\section{The Hyperbolic Double Doughnut.}

In this part we describe the construction of the Hyperbolic Double Doughnut. First we recall some important properties of the 2-dimensional
hyperbolic geometry. It is convenient to use the representation of the
hyperbolic space by using the Poincar\'e disc
\begin{equation}
\mathbf{D}:=\left\{(x,y)\in\RR^2,\;x^2+y^2<1\right\},
  \label{}
\end{equation}
endowed with the metric expressed with the polar coordinates by 
\begin{equation}
ds_{\mathbf{D}}^2=\frac{4}{(1-r^2)^2}\mbox{d}
r^2+4\frac{r^2}{(1-r^2)^2}\mbox{d}\varphi ^2=\frac{4}{(1-x^2-y^2)^2}
[\mbox{d}x^2+\mbox{d}y^2].
  \label{}
\end{equation}
It
is useful to use the complex parametrization $z=x+iy$. We have to
carefully distinguish the euclidean distance
\begin{equation}
d(z,z')=\mid z-z'\mid,
  \label{}
\end{equation}
and the hyperbolic distance associated with the hyperbolic metric,
given by
\begin{equation}
\cosh d_H(z,z')=1+\frac{2|z-z'|^2}{(1-|z|^2)(1-|z'|^2)}.
  \label{}
\end{equation}
We remark that
$$
d(0,z)=\tanh \frac{d_H(0,z)}{2}
$$
hence the euclidean circles centerd in $0$ are hyperbolic circles, and
more generally, all the hyperbolic circles $\{z';\;\;d_H(z',z_0)=R\}$,
with $R>0$, $z_0\in\mathbf{D}$, are euclidean circles.
The invariant measure $d\mu_H$ on the Poincar\'e disc allows to compute the
area of any Lebesgue measurable subset $X\subset\mathbf{D}$ by the formula
\begin{equation}
\mu_H(X)=\int_X \frac{4}{(1-|z|^2)^2}\; dx \; dy,
  \label{}
\end{equation}
in particular the hyperbolic area of a disc
$D_H(0,R):=\{z;\;d_H(z,0)\leq R\}$ is :
\begin{equation}
\mu_H(D_H(0,R))=4\pi \sinh^2\left(\frac{R}{2}\right).
  \label{}
\end{equation}

The group of the isometries of $\mathbf{D}$ is generated by three
kinds of transformations.

\begin{enumerate}

\item 
The Rotations of angle $\varphi_0\in\RR$ 
$$
R_{\varphi_0}(z)=e^{i\varphi_0} z,
$$
and so $R_{\varphi_0}$ is defined in the $(x,y)$-coordinates by the matrix
$$
R_{\varphi_0}=\left(\begin{array}{cc}
e^{i\frac{\varphi_0}{2}} & 0\\
0 & e^{-i\frac{\varphi_0}{2}}
\end{array}
\right).
$$

\item
The Boosts, or M\"obius transforms, associated with $\tau_0\in\RR$ : 
$$
T_{\tau_0}(z)=\frac{\cosh \frac{\tau_0}{2} z+\sinh
  \frac{\tau_0}{2}}{\sinh \frac{\tau_0}{2} z +\cosh \frac{\tau_0}{2}}
$$
expressed in $(x,y)$-coordinates by the matrix
$$
T_{\tau_0}=\left(\begin{array}{cc}
\cosh \frac{\tau_0}{2} & \sinh \frac{\tau_0}{2} \\
\sinh \frac{\tau_0}{2} &  \cosh \frac{\tau_0}{2}  
\end{array}
\right).
$$
We remark that 
$$
\cosh d_H
(z,T_{\tau_0}(z))=1+2\frac{|z^2-1|^2}{(1-|z|^2)^2}\sinh^2\left(\frac{\tau_0}{2}\right),
$$
therefore
\begin{equation}
\forall z\in]-1,1[,\;\;d_H (z,T_{\tau_0}(z))=\tau_0.
  \label{}
\end{equation}

\item
The Symmetry  
$$
S(z)=\overline {z}.
$$
\end{enumerate}

Finally we recall that the geodesics of the Poincar\'e disc are the diameters and all the arcs of
circles that intersect orthogonally the boundary of the disc.\\

Now we are ready to describe 
the double doughnut that is the quotient of the hyperbolic plane by the
Fuchsian group of isometries, $\Gamma$, generated by the four transforms
$g_0,g_1,g_2, g_3$, where
\begin{equation}
g_k\;=\;R_{k\frac{\pi}{4}}\;T_{\tau_1}\;R_{-k\frac{\pi}{4}}
  \label{}
\end{equation}
with
\begin{equation}
\tanh \frac{\tau_1}{2}=\sqrt {\sqrt{2}-1}.
  \label{}
\end{equation}
The matrix of $g_k$ is given by :
\begin{equation}
g_k=\left(
\begin{array}{cc}
1+\sqrt 2&\sqrt{2+2\sqrt 2}\,e^{ik\frac{\pi}{4}}\\
\sqrt{2+2\sqrt 2}e^{-ik\frac{\pi}{4}}&1+\sqrt 2
\end{array}
\right).
  \label{}
\end{equation}
These isometries $g_k$ satisfy the relation  : 
\begin{equation}
(g_0g_1^{-1}g_2g_3^{-1})(g_0^{-1}g_1g_2^{-1}g_3)\;=\;I_d.
  \label{}
\end{equation}
The {\it Hyperbolic Double Doughnut} is the quotient manifold
\begin{equation}
\mathbf{K}:=\mathbf{D}/\Gamma,
  \label{}
\end{equation}
endowed with the hyperbolic metric $ds_{\mathbf{K}}^2$ induced by
$ds_{\mathbf{D}}^2$ . We know, see e.g. \cite{Voroz}, \cite{bekka},
\cite{gutz},  that $\mathbf{K}$
is a two dimensional $C^{\infty}$ compact manifold, without boundary, its sectional
curvature is constant, equal to $-1$, and its genus, that is the number of ``holes'',
is 2. The geodesic flow is very chaotic : it is ergodic, mixing
(theorems by G.Hedlung, E. Hopf), Anosov and Bernouillian
(D. Ornstein, B. Weiss).\\

To perform the computations of the waves on the doughnut, it is very
useful to represent it by a minimal subset
$\mathcal{F}\subset\mathbf{D}$ and a relation of equivalence $\sim$
such that 
\begin{equation}
\mathbf{K}=\mathcal{F}/ \sim.
  \label{}
\end{equation}
 When $\mathcal{F}$ is
choosen as small as possible, it is called {\it Fundamental
  Polygon}. We take
\begin{equation}
\mathcal{F}:=\left\{ z\in\mathbf{D};\;\;\forall i=0,...,3, \; |g_{i}(z)|\geq|z|,\;\; |g_{i}^{-1}(z)|\geq|z|\right\}.
  \label{}
\end{equation}
We can see that $\mathcal{F} $ is a closed regular hyperbolic octogon,  of which the
boundary $\partial\mathcal{F}$
is the union of eight arcs of circle, that are parts of geodesics of $\mathbf{D}$. We denote $P_j$, $j\in\ZZ_8$, the
tops of  $\mathcal{F} $, and $\wideparen{P_jP_{j+1}}$ the eight wedges. The
action of $\Gamma$ on the boundary is described by :
\begin{equation}
\begin{array}{lllllll}
P_1 & = & g_3(P_6),  &  & P_5 & = & g_3^{-1}(P_2),\\
P_2 & = & g_2(P_7),  &  & P_6 & = & g_2^{-1}(P_3),\\
P_3 & = & g_1(P_8),  &  & P_7 & = & g_1^{-1}(P_4), \\
P_4 & = & g_0(P_1),  &  & P_8 & = & g_0^{-1}(P_5),
\end{array}
  \label{}
\end{equation}
and
\begin{equation}
\begin{array}{lllllll}
\wideparen{P_1 P_2} & = & g_3(\wideparen{P_6 P_5}), &  & \wideparen{P_2 P_3} & = & g_2(\wideparen{P_7 P_6}),\\
\wideparen{P_3 P_4} & = & g_1(\wideparen{P_8 P_7}), &  & \wideparen{P_4 P_5} & = & g_0(\wideparen{P_1 P_8}).
\end{array}
  \label{pp}
\end{equation}
\begin{figure}
\input{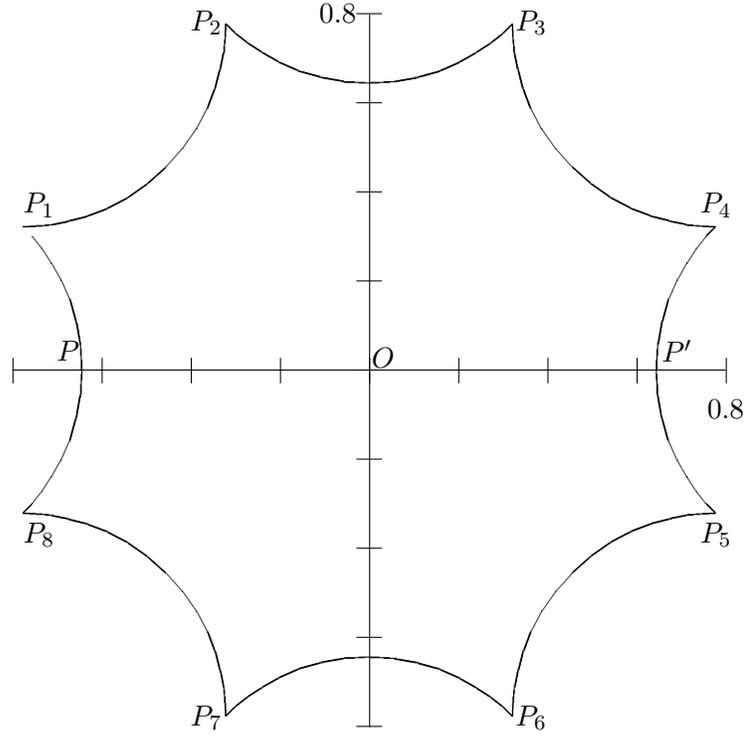}
\vspace{1cm}
\caption{The Fundamental Domain $\mathcal{F}$.}
\end{figure}

We define the relation $\sim$ by specifying the classes of equivalence
$\dot{z}$ of any $z\in\mathcal{F}$,
$\dot{z}:=\Gamma(\{z\})\cap\mathcal{F}$, i.e.
\begin{equation}
z \in \stackrel{\circ}{\mathcal{F}}\Rightarrow \dot{z}=\{z\},
  \label{}
\end{equation}
\begin{equation}
\dot{P_j}=\{P_1,\,P_2,\;P_3,\;P_4,\;P_5,\;P_6,\;P_7,\;P_8\},
  \label{}
\end{equation}
\begin{equation}
z\in\partial\mathcal{F}\setminus\dot{P_j}\Rightarrow \dot{z}=\{z,z_{equiv}\},
  \label{}
\end{equation}
where according (\ref{pp}) 
\begin{equation}
z\in\wideparen{P_iP_j},\;\;\wideparen{P_kP_l}=g_a^{\pm 1}\left(\wideparen{P_iP_j}\right)
\Rightarrow z_{equiv}=g_a^{\pm 1}(z)\in\wideparen{P_kP_l}.
  \label{equiv}
\end{equation}

We give some metric relations. We denote
$P=\wideparen{P_1P_8}\cap\RR^-$, $P'=g_0(P)=T_{\tau_1}(P)$. We have
$$
d_H(P,P')=2d_H(O,P')=\tau_1,
$$
and so
\begin{equation}
\wideparen{P_1P_8}\subset\left\{z;\;\left(x+\sqrt{\frac{1+\sqrt{2}}{2}}\right)^2 +y^2=\frac{\sqrt{2}-1}{2}\right\},
  \label{front}
\end{equation}
and by using the rotations we also have :
$$
d_H(P_i,P_{i+1})=\tau_1=2d_H(P,P_1)=2d_H(O,P),
$$
$$
d_H(O,P_i)=\tau_2\quad \mbox{with} \quad \tanh
\frac{\tau_2}{2}=2^{-\frac14},\;\;d(O,P_i)=2^{-\frac14}.
$$
Finally the area of the fundamental domain is $\mu_H(\mathcal{F})=4\pi$.


\section{The Waves on the Doughnut.}
The Laplace Beltrami operator  associated with a metric $g$, is defined by
$$
\frac{1}{\sqrt{|g|}} \partial_\mu g^{\mu \nu}
\sqrt{|g|}\partial_\nu,\;\;g^{-1}=(g^{\mu \nu})\;,\qquad |g|=|\det
g_{\mu \nu}|.
$$
We consider the lorentzian manifold
$\RR_t\times\mathbf{K}$ endowed with the metric 
\begin{equation}
g_{\mu\nu}dx^{\mu}dx^{\nu}=dt^2-ds_{\mathbf{K}}^2.
  \label{}
\end{equation}
We study the covariant wave equation
\begin{equation}
\partial_t^2\psi-\Delta_{\mathbf{K}}\psi=0,\;\;
  \label{}
\end{equation}
and more generally, the damped wave equation
\begin{equation}
\partial_t^2\psi-\Delta_{\mathbf{K}}\psi+a\partial_t\psi=0,\;\;
  \label{eq}
\end{equation}
whith $0\leq a\in L^{\infty}(\mathbf{K})$. Here $\Delta_{\mathbf{K}}$
  is the Laplace Beltrami operator on $\mathbf{K}$. Since $\mathbf{K}$ is a
smooth compact manifold without boundary, $\mathbf{K}$ endowed with
its natural domain $\{u\in L^2(\mathbf{K});\;\;\Delta_{\mathbf{K}}u\in
    L^2(\mathbf{K})\}$ is self-adjoint and the global Cauchy problem is
well posed in the framework of the finite energy spaces. Given
$m\in\NN$, we introduce the Sobolev space
\begin{equation}
  H^m(\mathbf{K}):=\left\{u\in
    L^2(\mathbf{K}),\;\nabla_H^{\alpha}u\in
    L^2(\mathbf{K}),\;\mid\alpha\mid\leq m\right\}
\end{equation}
where $\nabla_H$  are the covariant derivatives. We can also interpret
this space as the set of the distributions $u\in
H^m_{loc}(\mathbf{D})$ such that $u\circ g=u$ for any $g\in\Gamma$. Then for all
$\psi_0\in H^1(\mathbf{K})$, $\psi_1\in L^2(\mathbf{K})$, there exists
a unique $\psi\in C^0\left(\RR^+_t;H^1(\mathbf{K})\right)\cap
C^1\left(\RR^+_t;L^2(\mathbf{K})\right)$
solution of (\ref{eq}) satisfying
\begin{equation}
\psi(t=0)=\psi_0,\;\;\partial_t\psi(t=0)=\psi_1,
  \label{ci}
\end{equation}
and we have
\begin{equation}
\int_{\mathbf{K}}\mid\partial_t\psi(t)\mid^2+\mid\nabla_H\psi(t)\mid^2d\mu_H+\int_0^t\int_{\mathbf{K}}a\mid\partial_t\psi(t)\mid^2d\mu_Hdt=Cst.
  \label{}
\end{equation}

To perform the numerical computation of this solution,  we take
the fundamental polygon as the domain of calculus. Then the Cauchy problem
on $\RR_t\times\mathbf{K}$ is equivalent to the mixed problem
\begin{equation}
\partial_{tt}\psi-\frac{(1-x^2-y^2)^2}{4}\left[\frac{}{}\partial_{xx}\psi+\partial_{yy}\psi\right]+a(x,y)\partial_t\psi=0,\;\;(t,x,y)\in\RR^+\times\mathcal{F},
  \label{eqxy}
\end{equation}
with the boundary conditions
\begin{equation}
\forall (t,z)\in\RR\times\partial\mathcal{F},\;\;z\sim z'\Rightarrow \psi(t,z)=\psi(t,z').
  \label{cl}
\end{equation}
We denote $H^m(\mathcal{F})=\{u\in L^2(\mathcal{F}),\;\;\forall
\alpha\in\NN^2,\;\;\mid\alpha\mid\leq m,\;\;\partial_{x,y}^{\alpha}u\in L^2(\mathcal{F})\}$ the usual
Sobolev space $H^m$ for the euclidean metric, and we introduce
the spaces $W^m(\mathcal{F})$ that correspond to the spaces $H^m(\mathbf{K})$ :
\begin{equation}
W^m(\mathcal{F}):=\left\{u\vert_{\mathcal{F}};\;\;u\in H^m_{loc}(\mathbf{D}),\;\; \forall
  g\in\Gamma,\;\;u\circ g=u\right\},
  \label{}
\end{equation}
endowed with the norm
\begin{equation}
\|u\|_{W^m}:=\|u\vert_{\mathcal{F}}\|_{H^m(\mathcal{F})}.
  \label{}
\end{equation}
In particular, we have
\begin{equation}
W^1(\mathcal{F})=\left\{u\in H^1(\mathcal{F}),\;\;z\sim z'\Rightarrow
  u(z)=u(z')\right\},
  \label{}
\end{equation}
and for all
$\psi_0\in W^1(\mathcal{F})$, $\psi_1\in L^2(\mathcal{F})$, there exists
a unique $\psi\in C^0\left(\RR_t^+;W^1(\mathcal{F})\right)\cap
C^1\left(\RR_t^+;L^2(\mathcal{F})\right)$
solution of (\ref{ci}), (\ref{eqxy}) and (\ref{cl}), and we have the
energy estimate :
\begin{equation}
\begin{split}
 \int_{\F}\frac{4}{(1-x^2-y^2)^2}\mid\partial_{t}\psi(t,x,y)&\mid^2
 +\mid\partial_x\psi(t,x,y)\mid^2+\mid\partial_y\psi(t,x,y)\mid^2
 dxdy\\
&+\int_0^t \int_{\F}\frac{4a(x,y)}{(1-x^2-y^2)^2}\mid\partial_t\psi(t,x,y)\mid^2dxdy=Cst.
\end{split}
  \label{}
\end{equation}
Since $a\in L^{\infty}(\mathcal{F})$ we have a result of regularity : when $\psi_0\in W^2(\mathcal{F})$, $\psi_1\in W^1(\mathcal{F})$, then $\psi\in C^0\left(\RR_t;W^2(\mathcal{F})\right)\cap
C^1\left(\RR_t;W^1(\mathcal{F})\right)\cap
C^1\left(\RR_t;L^2(\mathcal{F})\right)$. In this case the mixed
problem can be expressed as a variational problem : $\psi$ is solution
iff for all $\phi\in W^1(\mathcal{F})$, we have :
\begin{equation*}
\begin{split}
\frac{d^2}{dt^2} \int_{\F}\frac{4}{(1-x^2-y^2)^2}\psi(t,z)\phi(z)&dxdy+
\frac{d}{dt} \int_{\F}\frac{4a(z)}{(1-x^2-y^2)^2}\psi(t,z)\phi(z)dxdy\\
&
-\int_{\F}\Delta_{x,y}\psi(t,z)\phi(z)dxdy=0.
\end{split}
  \label{}
\end{equation*}
To invoke the Green formula, we denote $\nu(z)$ the unit outgoing
normal at $z\in\partial\mathcal{F}$. We suppose that
$\wideparen{P_kP_l}=g_a\left(\wideparen{P_iP_j}\right)$. Then for $z\in\wideparen{P_iP_j}$
$$
g_a(\nu_z)=-\nu_{g_a(z)},
$$
Since $u\circ g_a=u$, we have for
$u\in W^2(\mathcal{F})$
$$
\partial_{\nu(z)}u(z)=g_a[\nu(z)].\nabla u(g_a(z))=-\partial_{\nu(g_a(z))}u(g_a(z)).
$$
We deduce that for $u\in W^2(\mathcal{F})$, $v\in W^1(\mathcal{F})$,
we have
$$\int_{\wideparen{P_iP_j}}v(z)\partial_{\nu(z)}u d\lambda(z)=
-\int_{\wideparen{P_kP_l}}v(z)\partial_{\nu(z)}u d\lambda(z),
$$
and therefore
$$
\int_{\partial\mathcal{F}}v(z)\partial_{\nu(z)}u d\lambda(z)=0,\;\;
\int_{\mathcal{F}}\Delta_{x,y}u(z)v(z)dxdy=-\int_{\mathcal{F}}\partial_xu\partial_xv+\partial_yu\partial_yvdxdy.
$$
We have proved the following

\begin{Theorem}
Given $\psi_0\in W^2(\mathcal{F})$, $\psi_1\in W^1(\mathcal{F})$, the
solution $\psi$ of the Cauchy problem (\ref{eq}), (\ref{ci}), is the
unique solution satisfying (\ref{ci}) of the variational problem
\begin{equation}
\begin{split}
\forall\phi\in W^1(\mathcal{F}),\;\;
\frac{d^2}{dt^2} \int_{\F}\frac{4}{(1-\mid z\mid^2)^2}\psi(t,z)\phi(z)&dxdy+
\frac{d}{dt} \int_{\F}\frac{4a(z)}{(1-\mid z\mid^2)^2}\psi(t,z)\phi(z)dxdy\\
&
+\int_{\F}\partial_x\psi(t,z)\partial_x\phi(z)+\partial_y\psi(t,z)\partial_y\phi(z)dxdy=0.
\end{split}
  \label{}
\end{equation}

  \label{}
\end{Theorem}

We solve this variational problem by the usual way. We take a family
$V_h$, $0<h\leq h_0$, of finite dimensional vector subspaces of
$W^1(\mathcal{F})$. We assume that
\begin{equation}
\overline{\cup_{0<h\leq h_0}V_h}=W^1(\mathcal{F}).
  \label{}
\end{equation}
We choose sequences $\psi_{0,h},\;\psi_{1,h}\in V_h$ such 
$$
\psi_{0,h}\rightarrow\psi_0\;\;in\;\;W^1(\mathcal{F}),\;\psi_{1,h}\rightarrow\psi_1\;\;in\;\;L^2(\mathcal{F}).
$$
We consider the solution $\psi_h\in C^{\infty}(\RR_t;V_h)$ of

\begin{equation}
\begin{split}
\forall\phi_h\in V_h,\;\;
\frac{d^2}{dt^2} \int_{\F}\frac{4}{(1-\mid z\mid^2)^2}\psi_h(t,z)\phi_h(z)&dxdy+
\frac{d}{dt} \int_{\F}\frac{4a(z)}{(1-\mid z\mid^2)^2}\psi_h(t,z)\phi_h(z)dxdy\\
&
+\int_{\F}\partial_x\psi_h(t,z)\partial_x\phi_h(z)+\partial_y\psi_h(t,z)\partial_y\phi_h(z)dxdy=0,
\end{split}
\end{equation}
satisfying $\psi_h(0,.)=\psi_{0,h}(.)$, $\partial_t\psi_h(0,.)=\psi_{1,h}(.)$.
Thanks to the conservation of the energy, this scheme is stable :
$$
\forall T>0,\;\;\sup_{0<h\leq h_0}\sup_{0\leq t\leq T}\|\psi_h(t)\|_{W^1}+\|\frac{d}{dt}\psi_h(t)\|_{L^2}<\infty.
$$
Moreover, when $\psi\in C^2\left(\RR^+_t;W^1(\mathcal{F})\right)$, it
is also converging :
$$
\forall T>0,\;\;\sup_{0\leq t\leq
  T}\|\psi_h(t)-\psi(t)\|_{W^1}+\|\frac{d}{dt}\psi_h(t)-\frac{d}{dt}\psi(t)\|_{L^2}\rightarrow 0,\;\;h\rightarrow 0.
$$
If we take a basis $\left(e_j^h\right)_{1\leq j\leq N_h}$ of $V_h$, we
expand $\psi_h$ on this basis :
$$
\psi_h(t)=\sum_{j=1}^{N_h}\psi_j^h(t)e_j^h,
$$
and we introduce 
$$
\MM=\left(M_{ij}\right)_{1\leq i,j\leq N_h},\;\;M_{ij}:=
 \int_{\F}\frac{4}{(1-\mid z\mid^2)^2}e_i^h(z)e_j^h(z)dxdy,
$$
$$
\DD=\left(D_{ij}\right)_{1\leq i,j\leq N_h},\;\;D_{ij}:=
 \int_{\F}\frac{4a(z)}{(1-\mid z\mid^2)^2}e_i^h(z)e_j^h(z)dxdy,
$$
$$
\KK=\left(K_{ij}\right)_{1\leq i,j\leq N_h},\;\;K_{ij}:=
 \int_{\F}\partial_xe_i^h(z)\partial_xe_j^h(z)+ \partial_ye_i^h(z)\partial_ye_j^h(z) dxdy,
$$
$$
X:=\left(
\begin{array}{c}
\psi_1^h\\
\psi_2^h\\
.\\
.\\
.\\
\psi_{N_h}^h
\end{array}
\right).
$$
Then the variational formulation is equivalent to 
\begin{equation}
\MM X''+\DD X'+\KK X=0.
  \label{}
\end{equation}
This differential system is solved very simply by iteration by solving
\begin{equation}
\MM(X^{n+1}-2X^n+X^{n-1})+\frac{\Delta
  T}{2}\DD(X^{n+1}-X^{n-1})+(\Delta T)^2\KK X^n=0.
  \label{}
\end{equation}
We know that this scheme is stable, and so convergent by the Lax
theorem, when
\begin{equation}
\sup_{X\neq 0}\frac{<\KK X,X>}{<\MM X,X>}< \frac{4}{\Delta T^2},
  \label{}
\end{equation}
and if there exists $K>0$ such that
\begin{equation}
\forall h\in]0,h_0],\;\;\forall \phi_h\in V_h,\;\;
\|\nabla_{x,y}\phi_h\|_{L^2(\mathcal{F})}\leq
\frac{K}{h}\left\|\frac{2}{1-\mid z\mid^2}\phi_h\right\|_{L^2(\mathcal{F})},
  \label{}
\end{equation}
the CFL condition
\begin{equation}
K\Delta T<2h,
  \label{}
\end{equation}
is sufficient to assure the stability and the convergence of our scheme.
\section{Numerical resolution}
\subsection{Mesh}

First of all we construct the boundary $\partial \mathcal{F} $ from
the equation (\ref{front}) and we perform a discretization that is equidistant for the hyperbolic metric.
Next we use the mesh generators Emc2\texttrademark and bamg\texttrademark created by INRIA.
If we only use Emc2\texttrademark, the mesh contains too many vertices and is not suitable for the hyperbolic metric.
 So a first mesh is created by Emc2\texttrademark. We also consider a circle which is uniformly discretized with the same hyperbolic step than the exterior geometry. The radius is choosen as the final mesh is almost uniform.
At last, we impose on every point of the exterior and interior
geometry a metric, in the sense of bamg\texttrademark. This software
can next create a mesh which is more uniform, with respect to the
hyperbolic metric, than the first mesh, and that  has a reasonnable number of vertices.

To test the uniformity of the mesh, we compute the extrema of the
hyperbolic distance between two neighbor vertices. As a check of the
accuracy of the meshes we evaluated the area of the polygons created
by the meshes, and we compared to $4 \pi$ (area of the domain). Here are some examples:

$$
\begin{array} {lcccl}
Mesh\; &number\ of\  vertices & \max d_H & \min d_H\ & area / 4 \pi \\
Mesh1 \;: &7448 &0.087&0.027&1.00012\\
Mesh2\;:& 17574 & 0.049 & 0.0177 &1.00007\\
Mesh3\;: &37329 &  0.036& 0.012& 1.00003\\
Mesh4\;: &67517  &0.027&0.009& 1.000018\\
 \end{array}
$$ 

In our meshes, the greater hyperbolic distance between consecutive
vertices is not reached near the exterior boundary. To give an idea of
the accuracy of this dicretization, we show in the following figure, a
very rough mesh : 

\begin{figure}[H]
   \centerline{\includegraphics[scale=0.5]{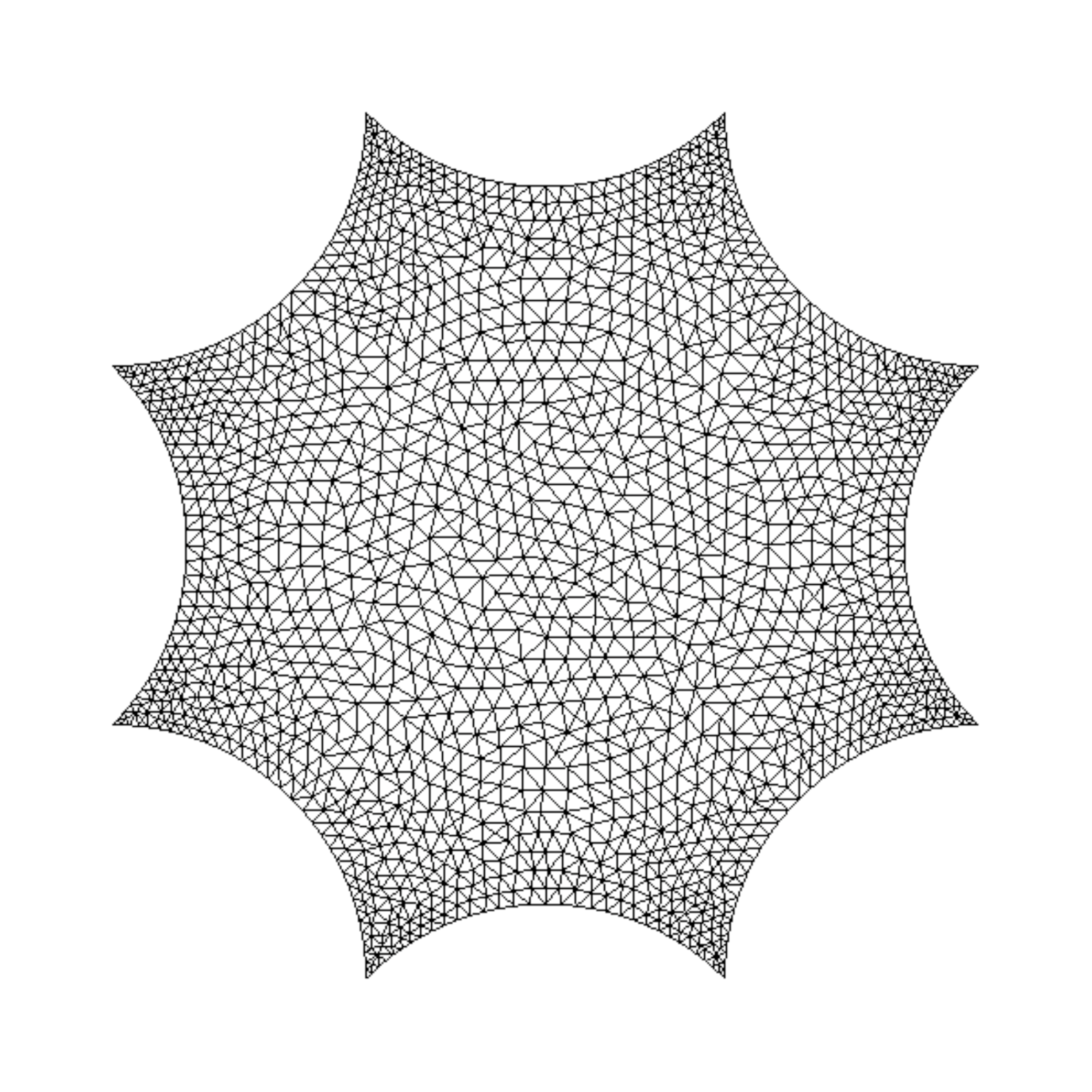}}
   \caption{A very rough mesh with 1756 vertices.}
\end{figure}

\subsection{$V_h$ space}
We construct the finite element spaces of $\P_1$ type. We note $\mathcal{T}_h$ all triangles of a mesh, and $\displaystyle{\mathcal{ F}_h=\cup_{K\in \mathcal{T}_h} K}$.
$$V_h:=\left\{ v: \mathcal{ F}_h \rightarrow \RR,\; v \in
\mathcal{C}^0(\mathcal{ F}_h),\; \forall K \in \mathcal{T}_h\;
v_{|K}\in \P_1(K),\; M\sim M'\Rightarrow v(M)=v(M')\right\}$$

If $M_i$ and $M_j$ denote two vertices of the mesh, we define a basis $\left(e_j^h\right)_{1\leq j\leq N_h}$ of $V_h$ by:\\
\begin{enumerate}
\item If $M_j\not \in \partial \mathcal{F}$ :
$e_j^h(M_i)=\delta _{ij}$ 
\item If $M_j$ is a $P_j$ point:
$e_j^h(M_i)=\left\{\begin{array}{ll}
1 & if\ M_i=P_i\\
0 & otherwise
\end{array}
\right.$
\item If $M_j\in \partial \mathcal{F}$, and is not a $P_j$ point: 
$e_j^h(M_i)=\left\{\begin{array}{ll}
1 & if\ M_i\sim M_j\\
0 & otherwise
\end{array}
\right.$

\end{enumerate}
In particular, we have to determine the equivalent points  on
$\partial \mathcal{F} $. To that, we write a program implementing the relations (\ref{equiv}).\\

The number of nodes $N_h$ is the sum of the number of the vertices which are not in $\partial \mathcal{F}$, the number of vertices which are on four consecutive arcs of $\partial \mathcal{F} $ without beeing a $P_i$ point, and one (because all $P_i$ points are equivalent to one of them).\\

\subsection{Matrix form of the problem}

$\KK(i,j)$ and $\MM(i,j)$ are found with a numerical integration using
the value at the middle of the edges of the triangles.\\

The stiffness matrix $\KK$ and the mass matrix $\MM$ are sparse and symetric matrices. So we choose a Morse stockage of their lower part, and all of the calculations will be performed with this stockage.\\

To solve the linear problem we use a preconditionned conjugate gradient method. The preconditionner is an incomplete Choleski factorisation, and the starting point is the solution obtained with a diagonal preconditionner. 

\subsection{Initial data}

We consider the case where the initial velocity $\psi_1=0$,
i.e. $X^0=X^1=0$, and we choose first initial datas with a more or
less small support near a given point. For instance, for  the wave depicted
in Figure 3, we have taken
\begin{equation}
\psi_0(x,y)=100e^{\frac{1}{100x^2+100y^2-1}},\;\;for\;x^2+y^2<
  \frac{1}{100},\;\psi_0(x,y)=0,\;\;for\;x^2+y^2\geq \frac{1}{100}.
  \label{data}
\end{equation}
\subsection{Discretized energy}

In order to see the stability of our method we perform $E_d(t)$ the discretized energy at the time $t$:

$$E(n\Delta t)=\left< \MM \;\frac{X^n-X^{n-1}}{\Delta t}\; , \;
  \frac{X^n-X^{n-1}}{\Delta t}\right>\;+\; \left< \KK\; X^{n-1}\; , \; X^n\right>$$
It is well known that our schema is conservative when $a=0$, hence $E_d$ must be
invariant all along the resolution. We test this property with the
previous initial data.

$$\begin{array} {|lc|cc|}
\hline
 &  E_d(0):& E_d (100): & \\
 & & &   \\
 &&\quad \rm{with\ }\Delta t=0.001: &\quad \rm{with\ }\ \Delta t=0.0005: \\
Mesh1 \;:&8455.89602005935&8455.89602005865&\\
Mesh2\;:&8484.17400988788 &8484.17400988815& 8484.17400988936    \\
Mesh3\;:&8494.43409419468&8494.43409419464& 8494.43409419511    \\
Mesh4 \;: &8498.39023175937  &8498.39023175945& 8498.39023175923 \\
\hline
 \end{array}
$$ 

\begin{figure}[H]
\includegraphics[scale=0.3]{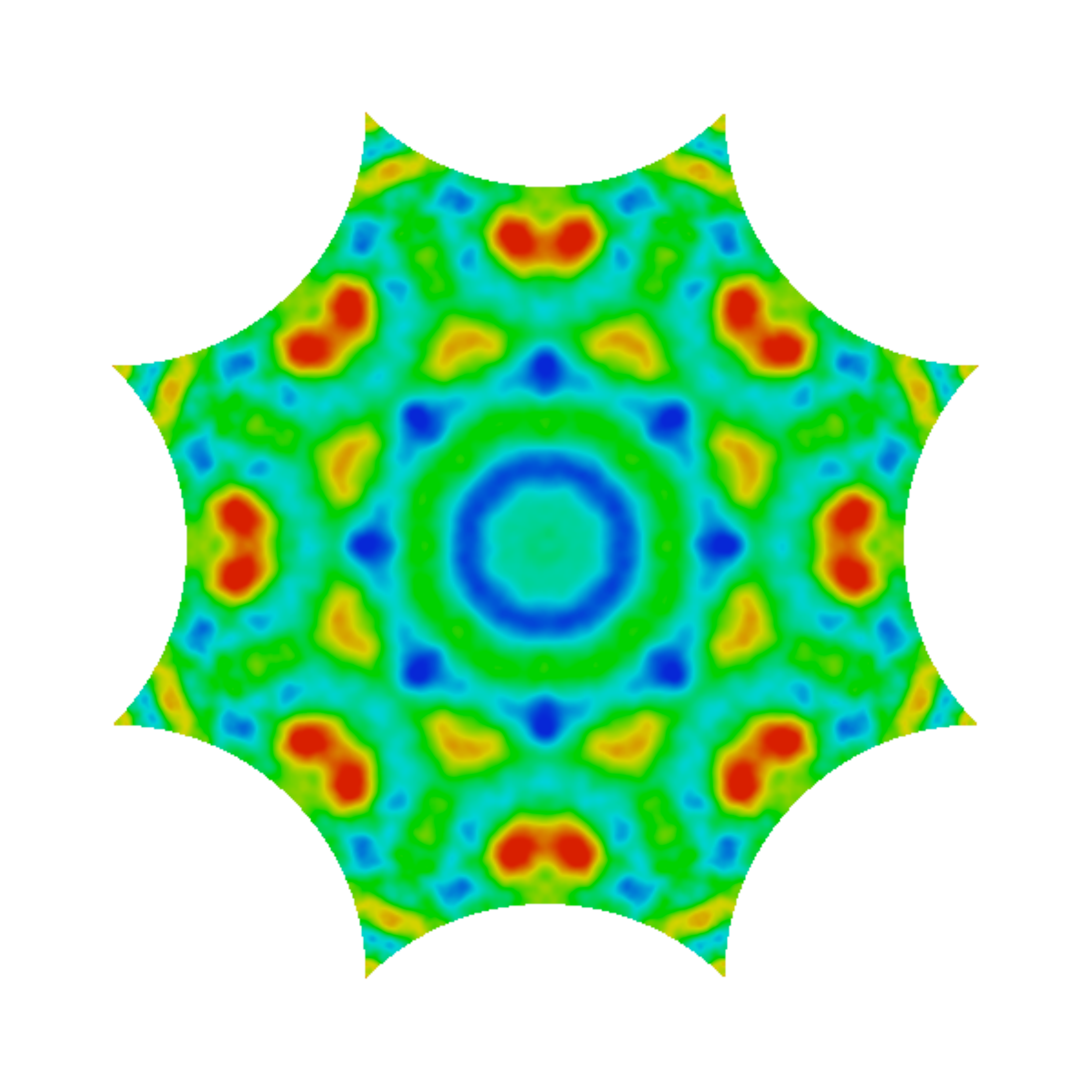}
\includegraphics[scale=0.3]{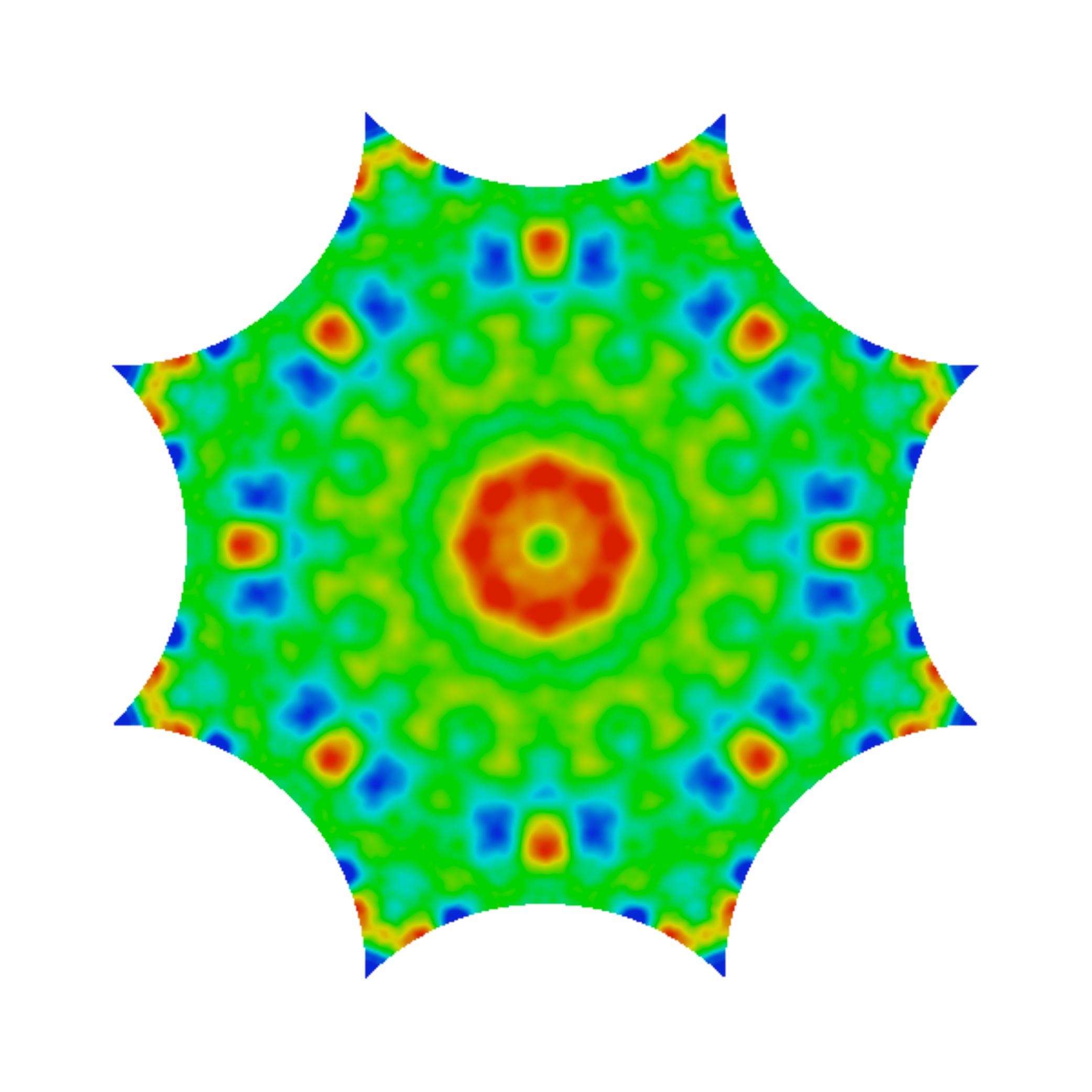}
\includegraphics[scale=0.3]{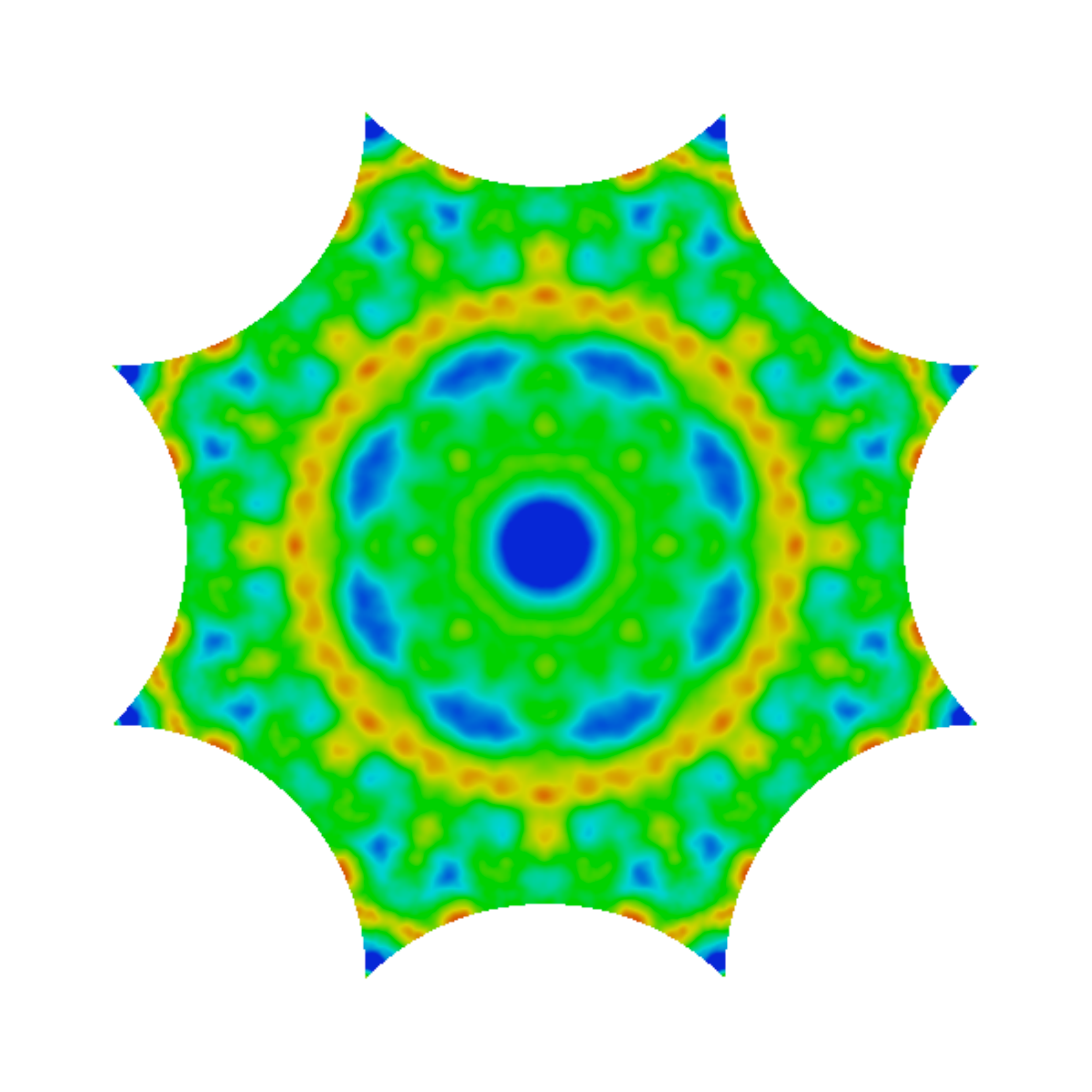}
\includegraphics[scale=0.3]{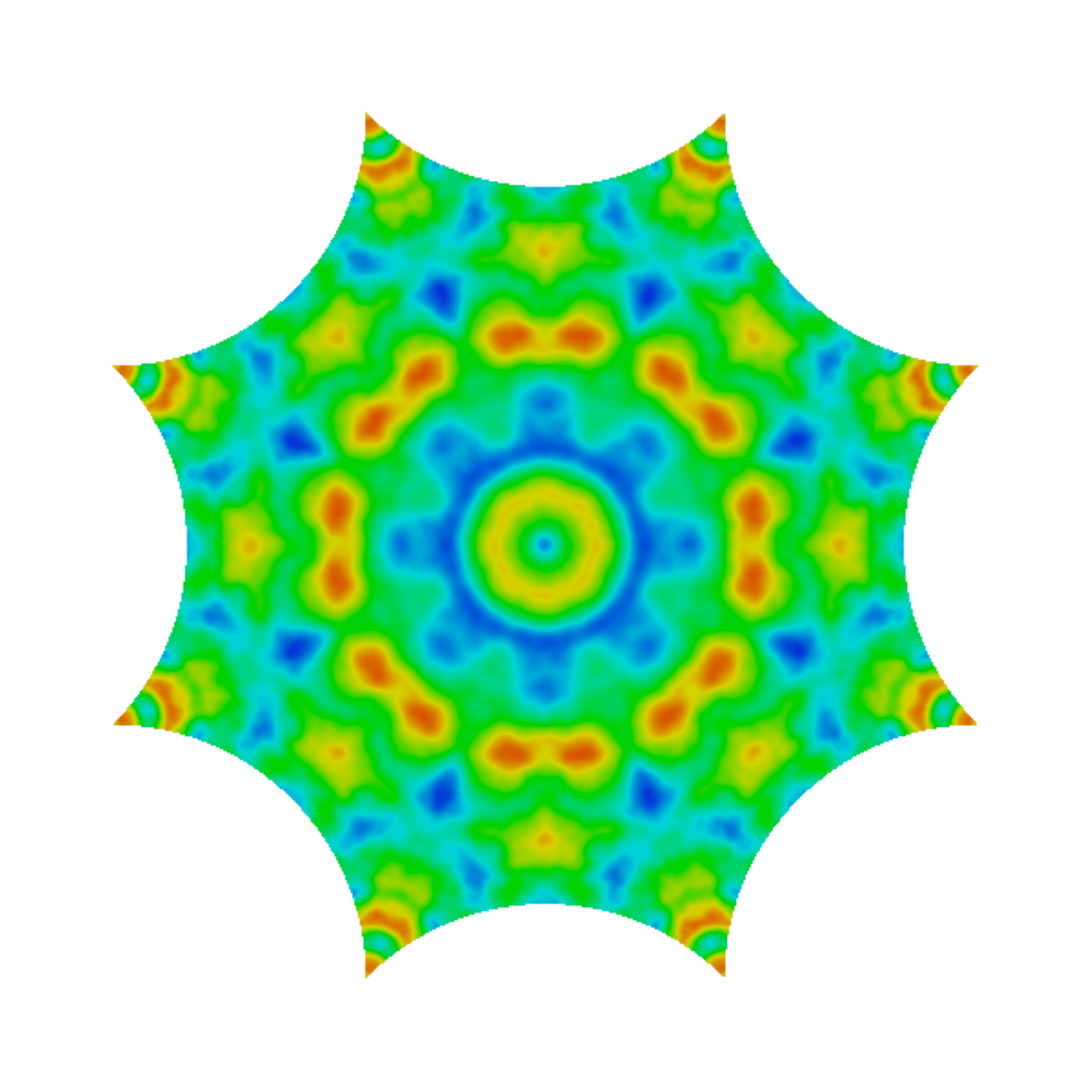}
\caption{The transient wave at $t=34$, $t=35$, $t=36$, $t=37$.}
\end{figure}


\subsection {Eigenvalues}
We test our scheme in the time domain, by looking for the eigenvalues
of the hamiltonian when $a=0$.
Since the Laplace-Beltrami operator $\Delta_{\mathbf{K}}$ on the
Hyperbolic Double Doughnut is a non positive, self-adjoint elliptic
operator on a compact manifold, its spectrum is a discrete set of
eigenvalues $-q^2\leq 0$,  and by the Hilbert-Schmidt theorem, there
exists an orthonormal basis in $L^2(\mathbf{K})$, formed of
eigenfunctions $\left(\psi_q\right)_q\subset H^{\infty}(\mathbf{K})$  associated to $q^2$, i.e.
\begin{equation}
-\frac{(1-x^2-y^2)^2}{4}\left[\frac{}{}\partial_{xx}\psi_q+\partial_{yy}\psi_q\right]=q^2\psi_q,\;\;\psi_q\in
W^{\infty}(\mathcal{F}).
  \label{}
\end{equation}
\noindent
We take $\psi_0=\frac{1}{2\sqrt\pi}$. Therefore any finite energy solution $\psi(t,x,y)$ of  $\partial_t^2
\psi-\Delta_{\mathbf{K}}\psi=0$ has an expansion of the  form $\sum_q
e^{iqt}\psi_q(x,y)$ (such expansions exist also for the damped wave
equation, when $a>0$, see \cite{hit}). 
More precisely, if we denote $<,>$ the scalar product in
$L^2(\mathbf{K})$, we write
\begin{equation}
\begin{split}
\psi(t,x,y)=&\frac{1}{4\pi}(<\partial_t\psi(0,.),1>t+<\psi(0,.),1>)\\
&+\sum_{q\neq 0} <\partial_t\psi(0,.),\psi_q>\frac{\sin
qt}{q}\; \psi_q(x,y)+\;<\psi(0,.),\psi_q>\cos qt\; \psi_q(x,y).
\end{split}
  \label{}
\end{equation}

To compute
the eigenvalues $q$, we investigate the Fourier transform in time of
the signal $\psi(t,x,y)$ in the case where $\partial_t\psi(0,x,y)=0$.
We fix some large $T>>1$, and we put
$\Psi_{\omega}(X):=\int_0^T \psi(t,X) e^{i\omega t} \; dt.$ Then
$\Psi_{\omega}(t)\sim C T^2(\omega^2-q^2)^{-1}$,
$T\rightarrow\infty$.
Practically, during the time resolution of the equation we store the
values of the solution at some points $M$, including the origin, $P'$,
$M_0$ near $P_4$, for the discrete time $k\Delta t$, $N_i\leq k\leq
N_f$. We choose the initial step $N_i$ in order to the transient wave
is stabilized, that to say $N_i\Delta t$ is greater than the diameter
of the doughnut, i.e. $N_i \Delta t\geq 2*d_H(0,P_i)\simeq 4.8969$.
Then we compute a DFT of $(\psi_h(k\Delta t,M))_{N_i\leq k \leq N_f}$
with the free FFT library fftw. Let us note
$(\Psi_j(M))_{0,N_f-N_i+1}$ the result. If $\Psi_j(M):=\sum_{k=0}^{N_f-N_i+1}\psi((N_i+k)\Delta t,x,y)e^{-(\frac{2\pi }{N_f-N_i+1}jk)i}$,
we search the values $jmax_1$, $jmax_2$, ... for which $(\|\Psi_j(M)\|^2)_j$ has a maximum. Then the eigenvalues found by the algoritm expressed as:
$$
q=\frac{2\pi}{(N_f-N_i+1)\Delta t}jmax
$$ 

We have made a lot of tests by varying parameters such as: $\Delta t$,
$N_i$, $N_f$, the mesh, the observation point
$M$. With the initial data (\ref{data}), we find the following
values for $q$ :

$$
\begin{array}{cc}
1,96\pm 0,02&
2,85\pm 0,04\\
4,34\pm 0,06&
4,83\pm 0,06\\
6,00\pm 0,06&
6,63\pm 0,06
\end{array}
$$

These results agree with the results obtained with a stationnary
method with a mesh of 3518 vertices in \cite{Aurich1}.\\
Alternatively, we could also use the power
spectrum and calculate the square of the modulus of $\Psi_{\omega}(x,y)$ 
\begin{equation*}
\begin{split}
\|\Psi_{\omega}\|_{L^2(\mathbf{K})}^2&=
\frac{1}{4\pi}\mid<\psi(0,.),1>\mid^2\\
&+\frac{1}{2}\sum_q <\psi(0,.),\psi_q>^2\left[\frac{\sin ^2
    \frac{(q+\omega)T}{2}}{(q+\omega)^2}+\frac{\sin ^2
    \frac{(\omega-q)T}{2}}{(\omega-q)^2}+\frac{1}{\omega
    ^2-q^2}\left(\cos ^2 \omega T-\cos \omega T \cos qT\right)\right]
\end{split}
  \label{}
\end{equation*}
Therefore, for an eigenvalue $q_0$: 
$$
 \|\Psi_{\omega}\|_{L^2(\mathbf{K})}^2
  \geq \; <\psi(0,.),\psi_q>^2\left( \frac{T^2}{4}-\frac{T}{2q}\right).
$$


\subsection{Damped waves}
We test our scheme for the damped wave equation (\ref{eq}) when the
damping function $a\geq 0$ is non zero (for deep theoretical results,
see \cite{hit}, \cite{lebeau}, \cite{Rauch}). We know that the energy
of any finite energy solution decays
exponentially (uniformly with respect to the initial energy) iff the dumping
$a$ satisfies the assumption of geometric control
introduced by J. Rauch and M. Taylor in
\cite{Rauch}. This condition means
\begin{equation}
\int_0^{\infty}a(x(t), y(t))dt=+\infty
  \label{cont}
\end{equation}
for any geodesic $(x(t), y(t))$. Since the geodesic flow on the
compact Riemaniann manifold with constant negative curvature is very
chaotic (more precisely ergodic, mixing, Anosov, Bernouillian see
e.g. \cite{Voroz}, \cite{bekka}, \cite{gutz}),
it is sufficient to have $a>0$ near $\partial\mathcal{F}$. Nevertheless
we constat an exponential decay for some solution, even if we choose a dumping function $a$
equal to a positive constant on
very small support that does not satisfy (\ref{cont}) : $a>0$ only on  one triangle and its close neighbors.

The next figures are obtained with mesh3 and $a$ defined by:
$$a(x,y)=0,\; {\rm for\ }, \mid z\mid<0.6 \quad;\quad a(x,y)=0.1,\;  {\rm otherwise.}$$
\begin{figure}[H]
\centering
 \input{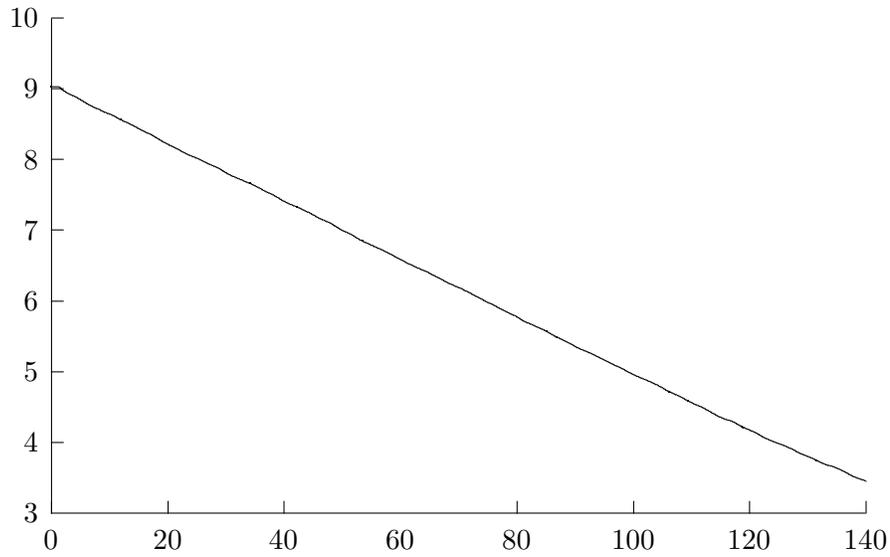}
\caption{Logarithm of the energy as a function of time.}
\end{figure}  
 
\begin{figure}[H]
\centering
 \input{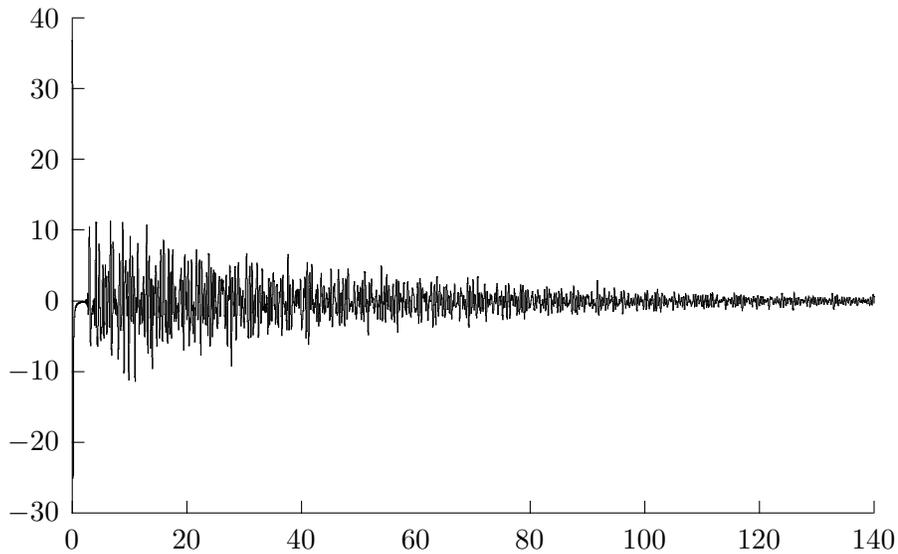}
\caption{Solution at the origin.}
\end{figure}


\begin{thebibliography}{5}

\bibitem{Voroz}
N.L. {\sc{Balazs}}, A. {\sc{Voros}}.
\newblock \protect {Chaos on the pseudosphere}.
\newblock {\em Phys. Rep.}, 143-3: 109-240, 1986.

\bibitem{bekka}
M. B. {\sc{Bekka}}, M. {\sc{Mayer}}.
\newblock \protect {Ergodic Theory and Topological Dynamics of Group
  Actions on Homogeneous Spaces}.
\newblock {\em London Mathematical Society Lecture Notes Series}, 269,
Cambridge University Press, 2000.

\bibitem{Cornish}
N.J. {\sc {Cornish}}, N.G. {\sc {Turok}}.
\newblock \protect {Ringing the eigenmodes from compact manifolds}.
\newblock {\em Class. Quantum Grav.}, 15:2699-2710, 1998.

\bibitem{Aurich1}
R. {\sc {Aurich}}, F. {\sc{Steiner}}.
\newblock \protect {Periodic-orbit sum rules for the Hadamard-Gutzwiller model}.
\newblock {\em Phys. D}, 39: 169-193, 1989.

\bibitem{gutz}
M. C. {\sc{Gutzwiller}}.
\newblock \protect {Chaos in Classical and Quantum Mechanics}.
\newblock {\em Interdisciplinary Applied Mathematics}, 1,
Springer-Verlag, 1990.


\bibitem{hit}
M. {\sc{Hitrik}}.
\newblock \protect {Eigenfrequencies and expansions for damped wave equations}.
\newblock {\em Methods Appl. Anal.}, 10,4 : 543-564, 2003.


\bibitem{Luminet}
M. {\sc {Lachi\`{e}ze-Rey}}, J.P. {\sc {Luminet}}.
\newblock \protect {}
\newblock {\em Phys. Rep.}, 254:135-214, 1995.

\bibitem{lebeau}
G. {\sc {Lebeau}}.
\newblock \protect {Equations des ondes amorties}.
\newblock {\em S\'eminaire X-EDP}, 15, Ecole Polytechnique, 1994.

\bibitem{Rauch}
J. {\sc {Rauch}}, M. {\sc {Taylor}}.
\newblock \protect {Decay of solutions to nondissipative hyperbolic systems on compact manifolds}.
\newblock {\em Comm. pure Appl. Math.}, 28:501-523, 1975.

\end{thebibliography}
\end{document}